# A Monolithic Eulerian Formulation for non-Classical Fluid-Structure Interaction (nCFSI): Modeling and Simulation

**Nazim Hussain[a], Muhammad Sabeel Khan[b,*], Lisheng Liu[a]**

[a] Department of Engineering Structure and Mechanics, Wuhan University of Technology, Wuhan, People's Republic of China, 430070.
[b] Department of Mathematics, Sukkur Institute of Business Administraiton University, 65200 Sukkur, Pakistan.

*Corresponding author*: m.sabeel@iba-suk.edu.pk

## Abstract
In this paper a new monolithic Eulerian formulation in the framework of non-classical continuum is presented for the analysis of fluid-strucutre interaction problems. In this respect, Cosserat continuum theory taking into account the micro-rotational degrees of freedom of the particles is considered. Continuum description of the model and variational formulation of the governing flow dynamics for non-classical -fluid-structure interaction nCFSI is presented. The model is analyzed by computing a well known benchmark problem by Hecht and Pironneau [16]. The algorithmic description is presented and implemented using FreeFEM++. Code is validated with the benchmark solution of Turek and Hron [38] in case of flow around a flag attached with cylinder. New microstructral behavior of the solution is studied and numerical simulations and results are shown in the form of figures. Some interesting feature of the flow is observed and microstructural characteristics are discussed.

## Graphical Abstract

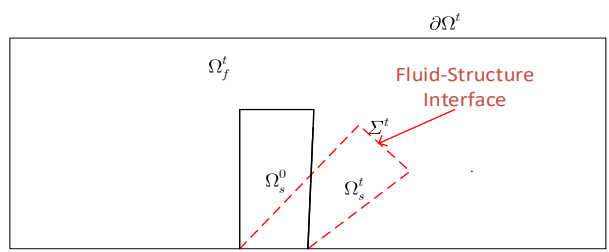
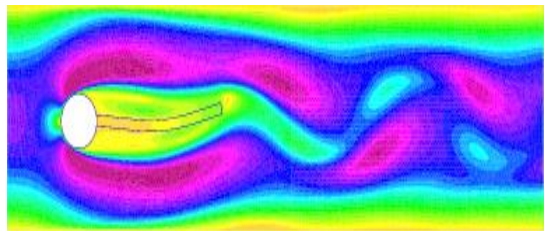

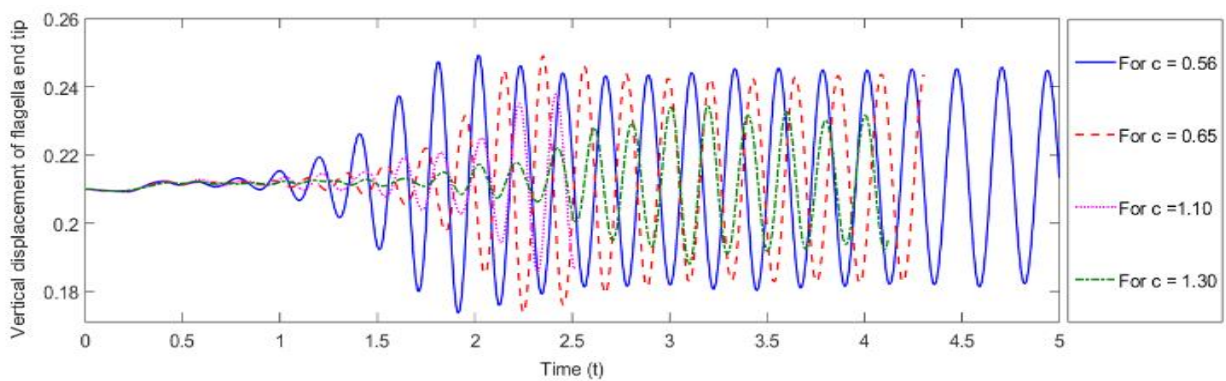

# 1 INTRODUCTION

Among the formulations used for solving fluid-structure interactions FSI problems some are based on monolithic approach, which finds origin in (Hron and Turek 2006) solve the FSI problems as a single variational equation for the whole system (Dunne 2006; Heil et al. 2008), or more recently (Pironneau 2016; Hecht and Pironneau 2017; Chiang et al. 2017; Pironneau 2018; Murea 2019). A monolithic Eulerian approach (Pironneau 2016; Hecht and Pironneau 2017) is similar to the fully Eulerian formulation (Dunne and Rannacher 2006; Dunne 2006), presents the governing equations of FSI in terms of displacement, velocity and pressure as usual but in an Eulerian framework. The major difference between two formulations is that in (Dunne and Rannacher 2006; Dunne 2006) prefer to work with velocities in the fluid domain and displacement in the solid domain, while monolithic Eulerian approach (Pironneau 2016; Hecht and Pironneau 2017) proposes to work with velocities everywhere in the problem domain. The fully Eulerian formulation replaces the well-established ALE, formulates both sub domains fluid and solid, respectively in Eulerian framework and well suited to problems with large displacements. More work devoted to this formulation can be seen in the previous studies such as (Richter and Wick 2010; Rannacher and Richter2011; Richter 2013; Wick 2013).

A partitioned approach is one in which the fluid and solid sub-problems are solved separately using iterative process: fixed point iterations (Formaggia et al. 2001; Nobile 2001; Tallec and Mouro 2001), Newton-like methods (Gerbeau and Vidrascu 2003; Fernández and Moubachir 2005; Dettmer and Perić 2006) or optimization techniques (Murea 2006; Mbaye and Murea 2008; Kuberry and Lee 2013).

Among the other formulations for analyzing FSI problems includes arbitrary Lagrangian-Eulerian formulation (ALE) approach used for simulation in the case of small displacement (Formaggia et al. 2010). The formulation matches the velocities and stresses at the fluid-structure interface, the fluid equations are then mapped back into solid domain at every time step during numerical simulation (Tallec and Hauret 2003; Hron and Turek 2006). However, in the case of large displacements the ALE formulation fails that leads to heavily distortion of fluid mesh (Liu 2016; Basting et al. 2017). Immersed boundary method IBM (Peskin 2002) is efficient for thin structures but implementation in thick structures remains challenging (Wang et al. 2017).

Numerical approaches exsisting in literature all are in classical continuum description. Here in this paper, a monolithic approach is presented but in the framework of non-classical Cosserat continuum description.

Classical continuum mechanics consider continuum as a simple point-continua with points having three displacement (dofs) and a symmetric Cauchy stress tensor characterizes the response of a material to the displacement. Such classical models may not be sufficient for the description of non-classical physical phenomena, where microstructural effects are observed most in high strain gradients regions. The Cosserat continuum (micropolar continuum) theory is one of the most prominent theory to model non-classical physical phenomena, and finds origin in (Eugene and François Cosserat 1909). Further, this concept was applied to describe fluids with microstructures by (Condiff and Dahler 1964; Eringen 1964 &1966) and the mathematical details, with some of its applications, are presented in (Lukaszewicz 1999). In non-classical continuum theory, the response of the material to the displacement and micro-rotation is characterized by a nonsymmetric Cauchy stress tensor and couple stresss tensor, respectively.

In the literature, the Cosserat fluids are not studied yet to analyze FSI phenomena due to high complexity. In the present study, a monolithic Eulerian approach is employed to analyze the Cosserat fluid-structure interaction CFSI. This benchmark problem was first studied by (Schafer and Turek 1996), and later by (Turek and Hron 2006; Dunne and Rannacher 2006; Hecht and Pironneau 2017), respectively. All the above considerations were in the case of Classical continuum but here in this paper a non-classical continuum framework is utilized to study the flow in particular the Cosserat continuum. The algorithmic description is presented and implemented using FreeFEM++ (Hecht 2012).

This paper is organized as follows. In section 2, we describe the Continuum description and notations used for the Cosserat model. Section 3, deals with the constitutive relations and the derivation of the governing equations from the conservation laws. In Section 4, we present the variational formulation in monolithic Eulerian framework for the CFSI problem. In Section 5, the time and spatial discretization is presented using semi-implicit scheme and the finite element method, respectively. Numerical tests and results are addressed in detailed in section 6. Future developments and conclusion of the study is addressed in section 7.

## 2 CONTINUUM DESCRIPTION AND NOTATIONS

Let $\Omega^t$ time dependent computational domain comprises of the fluid region $\Omega_f^t$ and solid region $\Omega_s^t$ such that $\overline{\Omega}^t = \overline{\Omega}_f^t \cup \overline{\Omega}_s^t, \Omega_f^t \cap \Omega_s^t = \varnothing, \forall t$. The interface of fluid and structure is denoted $\Sigma^t = \overline{\Omega}_f^t \cup \overline{\Omega}_s^t$ where the boundary of computational domain $\Omega^t$ is represented by $\partial \Omega^t$. At initial time the fluid domain $\Omega_f^0$ and the solid domain $\Omega_s^0$ are prescribed. The part of the boundary of computational domain $\partial \Omega^t$, where either the structure is clamped or the fluid has a no-slip condition is denoted by $\Gamma$ and it is assumed to be time independent. A schematic of this description is shown in Figure 1 below.

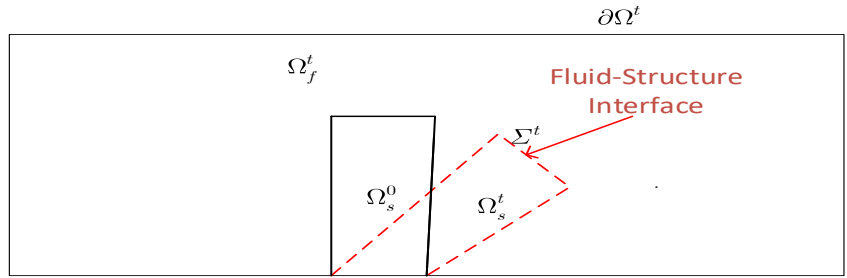

*Figure 1: Geometrical representation of the fluid-structure domain.*

Based on standard notations in (Ciarlet 1988; Marsden and Hughes 1993; Bath 1996; Tallec and Mouro 2001; Hron and Turek 2006), we consider $\mathbf{X}: \Omega^0 \times 0, T \mapsto \Omega^t : \mathbf{X}\ x^0,t$, the Lagrangian position at time $t$ of $x^0 \in \mathbb{R}^d, d = 2$ or $3$. The velocity of the deformation is given by $\mathbf{u} = \partial_t \mathbf{X}$ and $\mathbf{d} = \mathbf{X}\ x^0,t\ - x^0$ denotes the displacement vector of the structure. We represent the transposed gradient of deformation and the Jacobian of the deformation as $\mathbf{F}_{ji} = \partial_{x_i^0} \mathbf{X}_j$ and $J = \det \mathbf{F}$, respectively. Let the density and the stress tensor at given position $x$ and time $t$ are reprented as $\rho\ x,t$ and $\sigma\ x,t$. The density $\rho$ remains constant, in the fluid domain and in the solid domain for all time $t$, in the case of incompressible medium. Denoting the density constants by $\rho_f$ and $\rho_s$ at any point $x$ and time $t$, we define density by using the set function indicator $\mathbf{1}_{\Omega_f^t}$ as $\rho\ x,t\ = \rho_f \mathbf{1}_{\Omega_f^t}\ x,t\ + \rho_s \mathbf{1}_{\Omega_s^t}\ x,t$, where $\mathbf{1}_\Omega\ x\ = \begin{cases} 1 & if\ x \in \Omega \\ 0 & otherwise \end{cases}$.

Similarly the stress tensor at point $x$ and time $t$ is defined as $\boldsymbol{\sigma}\ x,t\ = \boldsymbol{\sigma}_f \mathbf{1}_{\Omega_f^t}\ x,t\ + \boldsymbol{\sigma}_s \mathbf{1}_{\Omega_s^t}\ x,t$. All spatial derivatives in this consideration are taken with respect to $x \in \Omega^t$ and not with respect to $x^0 \in \Omega^0$. If $\phi$ is a function of $x = \mathbf{X}\ x^0,t\$ where $x^0 \in \Omega^0$, then

$$\nabla_{x^0}\phi = \left[\partial_{x_i^0}\phi\right] = \left[\partial_{x_i^0}\mathbf{X}_j \partial_{x_j}\phi\right] = \mathbf{F}^T \nabla\phi. \tag{1}$$

The deformation gradient $\mathbf{F}$ and displacement field $\mathbf{d}$ can be seen as function of $x,t$ instead of $x^0,t$, and are related in the case when $\mathbf{X}$ is one-to-one and invertible by

$$\mathbf{F}^T = \nabla_{x^0}\mathbf{X} = \nabla_{x^0}(\mathbf{d}+x^0) = \nabla_{x^0}\mathbf{d}+\mathbf{I} = \mathbf{F}^T\nabla\mathbf{d}+\mathbf{I} \quad\Rightarrow\quad \mathbf{F} = (\mathbf{I}-\nabla\mathbf{d})^{-T}. \tag{2}$$

The time convective derivative of $\phi$ becomes

$$D_t\phi := \frac{\partial}{\partial t}\phi(\mathbf{X}(x^0,t),t) = \partial_t\phi(x,t) + \mathbf{u}\cdot\nabla\phi(x,t). \tag{3}$$

## 3 CONSTITUTIVE RELATIONS AND GOVERNING EQUATIONS

### 3.1 Consitutive Relation and Governing Equations for Solid medium

Conservation of mass and conservation of momentum for the solid medium takes form

$$\frac{d}{dt}(J\rho) = 0, \tag{4}$$

$$\rho D_t\mathbf{u} = \nabla\cdot\boldsymbol{\sigma}_s + \mathbf{f}. \tag{5}$$

Where $\mathbf{f}$ represents the density of volumetric forces and $J\rho = \rho^0$ in case of incompressible medium. For the structural part we assume an incompressible hyper-elastic material where constitutive description is stated as

$$\boldsymbol{\sigma}_s = -p_s\mathbf{I} + \rho_s\frac{\partial\Psi}{\partial\mathbf{F}}\mathbf{F}^T, \tag{6}$$

such that the Helmoltz potential $\Psi$ in the case of two dimensional Mooney-Rivilin material is defined as

$$\Psi(\mathbf{F}) = c_1 tr(\mathbf{F}^T\mathbf{F}) + c_2\left\{tr(\mathbf{F}^T\mathbf{F})^2 - tr^2(\mathbf{F}^T\mathbf{F})\right\}. \tag{7}$$

Where the constants $c_1$ and $c_2$ are empirically determined (Ciarlet 1988).

### 3.2 Derivation of Mooney-Rivlin 2-dimensional Stress Tensor

Since it is noted that $\partial_{\mathbf{F}}tr(\mathbf{F}^T\mathbf{F}) = 2\mathbf{F}$ and $\partial_{\mathbf{F}}tr(\mathbf{F}^T\mathbf{F})^2 = 4\mathbf{F}\mathbf{F}^T\mathbf{F}$. Hence

$$\partial_{\mathbf{F}}\Psi = 2c_1\mathbf{F} + c_2\left(4\mathbf{F}\mathbf{F}^T\mathbf{F} - 4tr(\mathbf{F}^T\mathbf{F})\mathbf{F}\right). \tag{8}$$

Let $\mathbf{B} = \mathbf{F}\mathbf{F}^T = \left[(\mathbf{I}-\nabla\mathbf{d})(\mathbf{I}-\nabla\mathbf{d})^T\right]^{-1}$, $b = \det\mathbf{B}$ and $c = tr\mathbf{B}$. Then

$$\partial_{\mathbf{F}}\Psi\mathbf{F}^T = (2c_1 - 4c_2 c)\mathbf{B} + 4c_2\mathbf{B}^2. \tag{9}$$

Now by applying the Caley-Hamilton theorem $\mathbf{B}^2 = c\mathbf{B} - b\mathbf{I}$,

$$\partial_{\mathbf{F}}\Psi\mathbf{F}^T = 2c_1\mathbf{B} - 4c_2 b\mathbf{I} = 2c_1\mathbf{F}\mathbf{F}^T - 4c_2\det(\mathbf{F}\mathbf{F}^T)\mathbf{I}. \tag{10}$$

Again by Caley-Hamilton theorem

$$\mathbf{B} = c\mathbf{I} - b\mathbf{B}^{-1} = c\mathbf{I} - b\left(\mathbf{I} - \nabla\mathbf{d} - \nabla\mathbf{d}^T + \nabla\mathbf{d}\nabla\mathbf{d}^T\right). \tag{11}$$

Using (11) in (10) one arrives at

$$\partial_{\mathbf{F}}\Psi\mathbf{F}^T = \left(2c_1(c-b)-4c_2 b\right)\mathbf{I}+2c_1 b\ \mathbf{Dd}-\nabla\mathbf{d}\nabla\mathbf{d}^T. \tag{12}$$

After some manipulations above equation becomes

$$\partial_{\mathbf{F}}\Psi\mathbf{F}^T = 2c_1 \det \mathbf{FF}^T \left(\mathbf{Dd}-\nabla\mathbf{d}\nabla\mathbf{d}^T\right) + \left(2c_1 tr\ \mathbf{FF}^T - 2c_1 + 4c_2\ \det \mathbf{FF}^T\right)\mathbf{I}. \tag{13}$$

Hence an incompressible $2$ - dimensional Mooney-Rivilin material will have

$$\partial_{\mathbf{F}}\Psi\mathbf{F}^T = 2c_1\left(\mathbf{Dd}-\nabla\mathbf{d}\nabla\mathbf{d}^T\right) + \alpha'\mathbf{I}. \tag{14}$$

Where $\alpha' = 2c_1 tr\ \mathbf{FF}^T - 2c_1 + 4c_2\ \det \mathbf{FF}^T$ is some scalar function of material parameters $c_1$ and $c_2$.

### 3.3 Consitutive Relation and Governing Equations for Cosserat Fluid

The mathematical formulation that governs the dynamics of Cosserat fluid is described as in the book of (Lukaszewicz 1999) and reads

$$\nabla \cdot \mathbf{u} = 0, \tag{15}$$

$$\rho\dot{\mathbf{u}} = \nabla \cdot \boldsymbol{\sigma}_f + \mathbf{f}, \tag{16}$$

$$\rho I\dot{\boldsymbol{\omega}} = \nabla \cdot \mathbf{C}_f + \boldsymbol{\varepsilon}:\boldsymbol{\sigma}_f + \mathbf{g}. \tag{17}$$

Where $\mathbf{u}$ and $\boldsymbol{\omega}$ are velocity and microrotation fields, respectively. The body force density, the body couple density and the microinertia coefficient are denoted by $\mathbf{f}$, $\mathbf{g}$ and $I$ respectively. We introduce the deformation tensor and the microrotaion strain tensor, respectively as

$$\mathbf{Du} = \nabla\mathbf{u}+\nabla\mathbf{u}^T, \tag{18}$$

$$\boldsymbol{\kappa} = \nabla\boldsymbol{\omega}. \tag{19}$$

For Cosserat fluid we considere the incompressible viscous relations described as

$$\boldsymbol{\sigma}_f = -p_f\mathbf{I}+\mu\left(\nabla\mathbf{u}+\nabla\mathbf{u}^T\right) +\mu_r\left(\nabla\mathbf{u}-\nabla\mathbf{u}^T\right) -2\mu_r\boldsymbol{\varepsilon}\cdot\boldsymbol{\omega}, \tag{20}$$

$$\mathbf{C}_f = \alpha\ tr\boldsymbol{\kappa}\ \mathbf{I}+\beta\left(\boldsymbol{\kappa}+\boldsymbol{\kappa}^T\right) +\gamma\left(\boldsymbol{\kappa}-\boldsymbol{\kappa}^T\right). \tag{21}$$

Where $\boldsymbol{\sigma}_f$, $\mathbf{C}_f$ and $\mathbf{I}$ are the nonsymmetric stress tensor, the couple stress tensor and the identity tensor, respectively. The pressure field and coefficient of dynamic viscosity are denoted by $p$ and $\mu$, respectively. The coefficients of microviscosity are represented as $\mu_r, \alpha, \beta$ and $\gamma$. The Levi Civita tensor is represented by $\varepsilon$. Subjected to certain prescribed boundary conditions according to the description of the physical problem and taking into consideration the constitutive equations (20) and (21), the governing conservation equations (15) - (17), leads to

$$\nabla \cdot \mathbf{u} = 0, \tag{22}$$

$$\rho\left(\frac{\partial \mathbf{u}}{\partial t} + \mathbf{u}\cdot\nabla\mathbf{u}\right) = -\nabla p + \mu + \mu_r \ \Delta\mathbf{u} + 2\mu_r \ \nabla\times\boldsymbol{\omega} + \mathbf{f}, \qquad (23)$$

$$\rho\left(I\frac{\partial \boldsymbol{\omega}}{\partial t} + I\mathbf{u}\cdot\nabla\boldsymbol{\omega}\right) = \lambda_1\Delta\boldsymbol{\omega} + \lambda_2\nabla\ \nabla\cdot\boldsymbol{\omega} - 4\mu_r\boldsymbol{\omega} + 2\mu_r(\nabla\times\mathbf{u}) + \mathbf{g}. \qquad (24)$$

Where $\lambda_1 = \beta + \gamma$ and $\lambda_2 = \alpha + \beta - \gamma$ are positive material parameter related to the microviscosity.

## 4 Monolithic Eulerian Variational Formulation of CFSI

In this variational formulation the homogeneous boundary conditions on $\Gamma \subset \partial\Omega^t$ are considere i.e., either solid is clamped or fluid has a no-slip condtition, and a homogeneous Neumann conditions on $\partial\Omega^t \setminus \Gamma$. In the case of incompressible material the Cosserat fluid-structure monolithic Eulerian variational formuation, thus reads

Given $\Omega_f^0$, $\Omega_s^0$ and $\mathbf{d}$, $\mathbf{u}$ at $t=0$, find $\mathbf{u}, \boldsymbol{\omega}, p, \mathbf{d}, \Omega_f^t, \Omega_s^t$ with $\mathbf{u}|_\Gamma = 0$ and $\boldsymbol{\omega}|_\Gamma = 0$, such that

$$\int_{\Omega_f \cup \Omega_s} \rho D_t\mathbf{u}\cdot\tilde{\mathbf{u}} - p\nabla\cdot\tilde{\mathbf{u}} - \tilde{p}\nabla\cdot\mathbf{u} + (\mu+\mu_r)\mathbf{Du}:\mathbf{D\tilde{u}} - 2\mu_r(\nabla\times\boldsymbol{\omega})\cdot\tilde{\mathbf{u}}\ d\Omega$$
$$+ \int_{\Omega_s} c_3(\mathbf{Dd} - \nabla\mathbf{d}\nabla\mathbf{d}^T):\mathbf{D\tilde{u}}\,d\Omega_s = \int_{\Omega_f \cup \Omega_s} \mathbf{f}\cdot\tilde{\mathbf{u}}\,d\Omega, \qquad (25)$$

$$\int_{\Omega_f \cup \Omega_s} \rho D_t\ I\boldsymbol{\omega}\cdot\tilde{\boldsymbol{\omega}} + \lambda_1\ \nabla\boldsymbol{\omega}:\nabla\tilde{\boldsymbol{\omega}} - \lambda_2\nabla\cdot\boldsymbol{\omega}\cdot\tilde{\boldsymbol{\omega}} + 4\mu_r\boldsymbol{\omega}\cdot\tilde{\boldsymbol{\omega}} - 2\mu_r\ \nabla\times\mathbf{u}\cdot\tilde{\boldsymbol{\omega}}\ d\Omega = \int_{\Omega_f \cup \Omega_s} \mathbf{g}\cdot\tilde{\boldsymbol{\omega}}\,d\Omega, \qquad (26)$$

$\forall(\tilde{\mathbf{u}}, \tilde{\boldsymbol{\omega}}, \tilde{p})$ taking $\tilde{\mathbf{u}}|_\Gamma = 0$ and $\tilde{\boldsymbol{\omega}}|_\Gamma = 0$, where $\Omega_f^t$ and $\Omega_s^t$ are defined incrementally by, $D_t\mathbf{d} = \mathbf{u}$ and

$$\frac{dX}{d\tau} = \mathbf{u}(X(\tau),\tau), \quad X(t)\in\Omega_r^t \Rightarrow X(\tau)\in\Omega_r^\tau \quad \forall\tau\in(0,T), \quad r=s,f. \qquad (27)$$

The relation in (27) defines $\Omega_f^t$ and $\Omega_s^t$ forward in time. Above the notations $\mathbf{B}:\mathbf{C} = tr(\mathbf{B}^T\mathbf{C})$ and $c_3 = \rho^s c_1$ are used.

## 5 DISCRETIZATION

In this section, we present the discretization scheme used to approximate the CFSI problem in equations (25-27). We use a semi-implicit scheme for time discretization and Galerkin discretization finite element method for space.

### 5.1 A Monolithic Time – Discrete Formulation of CFSI

Let $t\in[0,T]$ be the time of simulation where $T$ is the total time. Discretize the interval $[0,T]$ into equall sub intervals each of the length $\delta t = \frac{T}{N}$ such that $t = n\delta t$ where $n = 0,1,\cdots,N$. Let $\mathbf{d}^{n+1} = \mathbf{d}^n + \delta t\mathbf{u}^{n+1}$. Hence

$$\mathbf{Dd} - \nabla\mathbf{d}\nabla\mathbf{d}^T \approx \mathbf{Dd}^n - \nabla\mathbf{d}^n{\mathbf{d}^n}^T + \delta t\ \mathbf{Du}^{n+1} - \nabla\mathbf{u}^{n+1}{\mathbf{d}^n}^T - \nabla\mathbf{d}^n{\nabla\mathbf{u}^{n+1}}^T + o\ \delta t\ . \qquad (28)$$

Now, if $X^n$ is a first order approximation of $X(t^{n+1} - \delta t)$ defined by $\dot{X} = \mathbf{u}(X(\tau),\tau)$, $X(t^{n+1}) = x$, where $X(t^{n+1}) = x$ such that $X^n(x) = x - \delta t \mathbf{u}^n(x)$, then a first order in time approximation for CFSI problem (25-26) reads:

Find $\mathbf{u}^{n+1} \in H_0^1(\Omega^{n+1})$, $\boldsymbol{\omega}^{n+1} \in H_0^1(\Omega^{n+1})$, $p^{n+1} \in L_0^2(\Omega^{n+1})$, $\Omega_f^{n+1}$ and $\Omega_s^{n+1}$ such that with $\mathbf{u}^{n+1}|_\Gamma = 0$, $\boldsymbol{\omega}^{n+1}|_\Gamma = 0$ and $\Omega^{n+1} = \Omega_f^{n+1} \cup \Omega_s^{n+1}$, $\forall \, \tilde{\mathbf{u}} \in H_0^1(\Omega^{n+1}), \tilde{\boldsymbol{\omega}} \in H_0^1(\Omega^{n+1}), \tilde{p} \in L_0^2(\Omega^{n+1})$ with $\tilde{\mathbf{u}}|_\Gamma = 0$ and $\tilde{\boldsymbol{\omega}}|_\Gamma = 0$, the following holds

$$\int_{\Omega_f^n \cup \Omega_s^n} \left\{ \left( \rho^n \frac{\mathbf{u}^{n+1} - \mathbf{u}^n \circ X^n}{\delta t} \right) \cdot \tilde{\mathbf{u}} - p^{n+1} \nabla \cdot \tilde{\mathbf{u}} - \tilde{p} \nabla \cdot \mathbf{u}^{n+1} + (\mu + \mu_r) \mathbf{D}\mathbf{u}^{n+1} : \mathbf{D}\tilde{\mathbf{u}} - 2\mu_r (\nabla \times \boldsymbol{\omega}^{n+1}) \cdot \tilde{\mathbf{u}} \right\} d\Omega^n$$
$$+ \int_{\Omega_s^n} c_3 \left[ \mathbf{D}\mathbf{d}^n - \nabla \mathbf{d}^{n^T} \nabla \mathbf{d}^n + \delta t \left( \mathbf{D}\mathbf{u}^{n+1} - \nabla \mathbf{u}^{n+1} \nabla \mathbf{d}^{n^T} - \nabla \mathbf{d}^n \nabla \mathbf{u}^{n+1^T} \right) : \mathbf{D}\tilde{\mathbf{u}} \right] d\Omega_s^n = \int_{\Omega_f^n \cup \Omega_s^n} \mathbf{f} \cdot \tilde{\mathbf{u}} \, d\Omega^n, \quad (29)$$

and

$$\int_{\Omega_f^n \cup \Omega_s^n} \left\{ \rho^n I^n \left( \frac{\boldsymbol{\omega}^{n+1} - \boldsymbol{\omega}^n \circ \tilde{X}^n}{\delta t} \right) \cdot \tilde{\boldsymbol{\omega}} + \lambda_1 \nabla \boldsymbol{\omega}^{n+1} : \nabla \tilde{\boldsymbol{\omega}} - \lambda_2 \nabla (\nabla \cdot \boldsymbol{\omega}^{n+1}) \cdot \tilde{\boldsymbol{\omega}} + 4\mu_r \boldsymbol{\omega}^{n+1} \cdot \tilde{\boldsymbol{\omega}} - 2\mu_r (\nabla \times \mathbf{u}^{n+1}) \cdot \tilde{\boldsymbol{\omega}} \right\} d\Omega^n$$
$$= \int_{\Omega_f^n \cup \Omega_s^n} \mathbf{g} \cdot \tilde{\boldsymbol{\omega}} \, d\Omega^n. \quad (30)$$

Now, update $\mathbf{d}$ by $\mathbf{d}^{n+1} = \mathbf{d}^n \circ X^n + \delta t \mathbf{u}^{n+1}$, and $\Omega_r^n$ by $\Omega_r^{n+1} = \{x + \delta t \mathbf{u}^{n+1}(x) : x \in \Omega_r^n\}$, where $r = s, f$.

### 5.2 A Monolithic Spatial – Discretization with Finite Elements of CFSI

Let $V_h$ and $W_h$ represents the finite element functional spaces for the velocities, displacements and micro-rotational velocities, respectively and $Q_h$ be the functional space for pressure field. Let $\mathfrak{S}_h^0$ be a triangulation of the initial domain $\Omega^0$ with quadratic elements for displacements, translational, micro-rotational velocities and linear elements for pressure field. Given that that the pressure is different in fluid domain and structural domain because of the discontinuity of pressure at the fluid-structure interface $\Sigma$; therefore, the functional space $Q_h$ is space of piecewise linear functions on the triangulation and is continuous in $\Omega_r^{n+1}, r = s, f$. A small penalization parameter $\zeta \ll 1$ needs to be added to impose uniqueness of the pressure when one desire to use direct linear solver. The discrete variational formulation of CFSI, thus reads

Find $\mathbf{u}_h^{n+1}, \boldsymbol{\omega}_h^{n+1}, p_h^{n+1} : \forall \, \tilde{\mathbf{u}}_h \in V_{0h}, \tilde{\boldsymbol{\omega}}_h \in \tilde{W}_{0h}, \tilde{p}_h \in Q_h$ with $V_{0h}|_\Gamma = 0$ and $W_{0h}|_\Gamma = 0$ are subspaces of $V_h$ and $W_h$, such that

$$\int_{\Omega_f^n \cup \Omega_s^n} \left\{ \rho^n \frac{\mathbf{u}_h^{n+1} - \mathbf{u}_h^n \circ X^n}{\delta t} \cdot \tilde{\mathbf{u}}_h - p_h^{n+1} \nabla \cdot \tilde{\mathbf{u}}_h - \tilde{p}_h \nabla \cdot \mathbf{u}_h^{n+1} + (\mu + \mu_r) \mathbf{D}\mathbf{u}_h^{n+1} : \mathbf{D}\tilde{\mathbf{u}}_h - 2\mu_r (\nabla \times \boldsymbol{\omega}_h^{n+1}) \cdot \tilde{\mathbf{u}}_h \right\} \Omega^n$$
$$+ \int_{\Omega_s^n} c_3 \left[ \mathbf{D}\mathbf{d}_h^n - \nabla \mathbf{d}_h^{n^T} \nabla \mathbf{d}_h^n + \delta t \left( \mathbf{D}\mathbf{u}_h^{n+1} - \nabla \mathbf{u}_h^{n+1} \nabla \mathbf{d}_h^{n^T} - \nabla \mathbf{d}_h^n \nabla \mathbf{u}_h^{n+1^T} \right) : \mathbf{D}\tilde{\mathbf{u}}_h \right] d\Omega_s^n + \int_{\Omega_f^n \cup \Omega_s^n} \zeta p_h \tilde{p}_h \, d\Omega^n = \int_{\Omega_f^n \cup \Omega_s^n} \mathbf{f} \cdot \tilde{\mathbf{u}}_h \, d\Omega^n, \quad (31)$$

and

$$\int_{\Omega_f^n \cup \Omega_s^n} \left\{ \rho^n I^n \left( \frac{\boldsymbol{\omega}_h^{n+1} - \boldsymbol{\omega}_h^n \circ \tilde{X}^n}{\delta t} \right) \cdot \tilde{\boldsymbol{\omega}}_h + \lambda_1 \nabla \boldsymbol{\omega}_h^{n+1} : \nabla \tilde{\boldsymbol{\omega}}_h - \lambda_2 \nabla \left( \nabla \cdot \boldsymbol{\omega}_h^{n+1} \right) \cdot \tilde{\boldsymbol{\omega}}_h + 4\mu_r \boldsymbol{\omega}_h^{n+1} \cdot \tilde{\boldsymbol{\omega}}_h - 2\mu_r \left( \nabla \times \mathbf{u}_h^{n+1} \right) \cdot \tilde{\boldsymbol{\omega}}_h \right\} d\Omega^n \qquad (32)$$

$$= \int_{\Omega_f^n \cup \Omega_s^n} \mathbf{g} \cdot \tilde{\boldsymbol{\omega}}_h d\Omega^n.$$

To update the triangulation at each vertex (say $q_i^n$) of the triangle $T_h \in \mathfrak{I}_h^n$, the vertex is moved to a new position by

$$q_i^{n+1} := q_i^n + \delta t \mathbf{u}_h^{n+1}.$$

By denoting $\mathbf{d}_i^n := \mathbf{d}^n(q_i)$, it can be seen that

$$\mathbf{d}^n \circ X^n (q_i^{n+1}) = \mathbf{d}^n (q_i^n + \delta t \mathbf{u}_h^{n+1} - \delta t \mathbf{u}_h^{n+1}) = \mathbf{d}^n (q_i^n).$$

This implies that the displacement vector of vertices $\mathbf{d}_h^n$ can be copied to $\mathbf{d}_i^{n+1}$ plus with addition of $\delta t \mathbf{u}_h^{n+1}(q_i^n)$ in order to obtain $\mathbf{d}_h^{n+1}$, i.e.,

$$\mathbf{d}_h^{n+1} = \mathbf{d}_h^n \circ X^n + \delta t \mathbf{u}_h^{n+1} (q_i^n) = \mathbf{d}_h^n + \delta t \mathbf{u}_h^{n+1} (q_i^n).$$

Moreover, the fluid domain mesh is moved by $\tilde{\mathbf{u}}$ which is a solution of the Laplace problem $-\Delta \tilde{\mathbf{u}} = 0$ $\forall \tilde{\mathbf{u}} \in V_{0h}$, subjected to $\tilde{\mathbf{u}}|_\Sigma = \mathbf{u}$ where $\Sigma$ is Cosserat fluid structure interface and $\tilde{\mathbf{u}} = 0$ at the boundaries $\Gamma^f \cup \Gamma^s \setminus \Sigma$. Moving the vertices of each triangle $T_h \in \mathfrak{I}_h^n$ by the above procedure gives a new triangulation $\mathfrak{I}_h^{n+1}$.

## 6 RESULTS AND DISCUSSION

The monolithic Eulerian formulation is used to compare the non-classical case of Cosserat fluid-structure interaction CFSI problem with the classical case of fluid-structure interaction FSI. The space discretization is done by using Lagrangian triangular finite elements with quadratic elements for displacements, translational, micro-rotational velocities and linear elements for pressure field. The public domain software FreeFEM++ (Hecht 2012) has been used to implement the algorithms. The presented model is analyzed by computing a well known benchmark problem FLUSTUK-FSI-3* (Hecht and Pironneau 2017). This benchmark problem was first studied by (Schafer and Turek 1996), and later by (Turek and Hron 2006; Dunne and Rannacher 2006; Hecht and Pironneau 2017), respectively. The description of the model problem in consideration is shown schematically in Figure 2. An incompressible hyperelastic Mooney-Rivlin material, like rectangular flag of size $[0,l] \times [0,h]$, is attached at the back of a hard fixed cylinder in the computational rectangular domain $[0,L] \times [0,H]$; the fluid flow enters and leaves freely at the left inlet and the right outlet, respectively. The configuration, boundary conditions and intial conditions for the present test problem are described as follows.

**Configuration**

The point $0.2, 0.2$ is the center of a cylinder; other parameters are considered in computation as $l = 0.35$, $h = 0.02$, $L = 2.5$ and $H = 0.41$ which set cylinder slightly below the symmetry line of the computational domain.

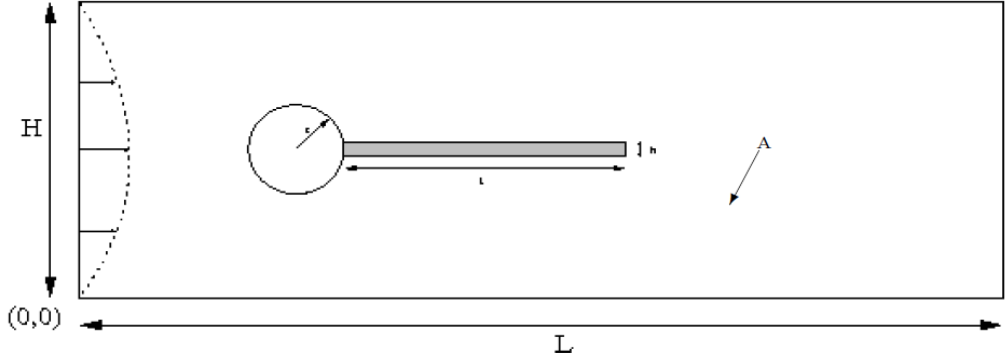

*Figure 2: Computational domain of the model problem.*

**Boundary and Initial Conditions**

Top and bottom boundaries satisfies the 'no-slip' condition. A parabolic velocity profile is prescribed at the left inlet,

$$\mathbf{u}_f(0, y) = \bar{U}\left(\frac{6y(H-y)}{H^2}\right),$$

where $\bar{U}$ is mean inflow velocity with flux $\bar{U}H$ and $\bar{U} = 2$. The zero-stress $\sigma \cdot n = 0$ is employed at the right outlet using do-nothing approach. Furthermore, the density and the reduced kinematic viscosity of the fluid takes values $\rho_f = 10^3 \, kgm^{-3}$ and $\nu_f = \frac{\mu}{\rho_f} = 10^{-3} m^2 s^{-1}$. For solid structure we consider $\rho_s = \rho_f$, $c_1 = 10^6 \, kgm^{-1}s^{-1}$ and no external force. Initially, all flow velocities and structure displacements are zero.

The flow starts oscillating and develops a Karman vortex street around time $t \sim 2$. The Figure 3 and Figure 4, shows the horizontal and vertical displacement of flagella end tip as a function of time in non-classical sense for CFSI. The results are obtained using mesh of 2199 vertices and time step size to 0.005 as in (Hecht and Pironneau 2017). The obtained results are compared and validated with the (Hecht and Pironneau 2017), showing amplitude and frequency roughly same.

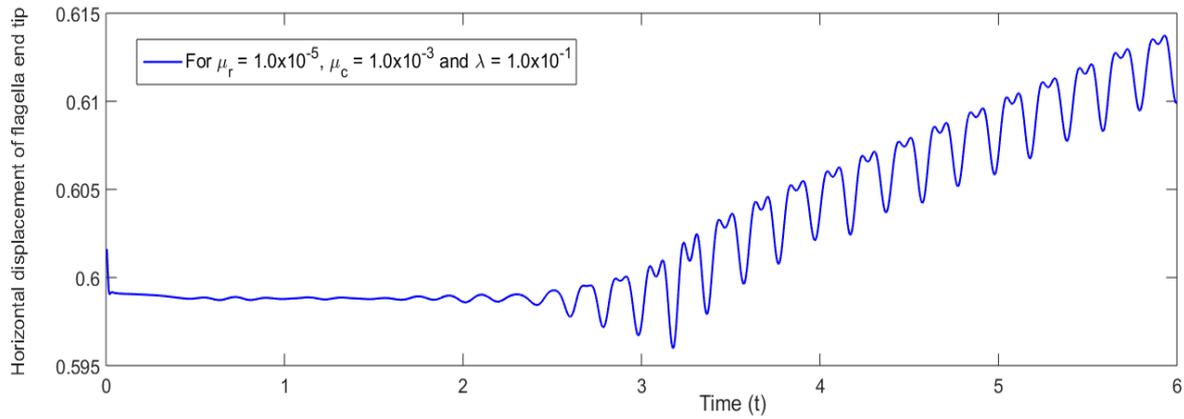

**Figure 3**: *Horizontal displacement of flagella end tip (control point A) against microrotional viscosity $\mu_r$, microinnertia $\mu_c$ at different time.*

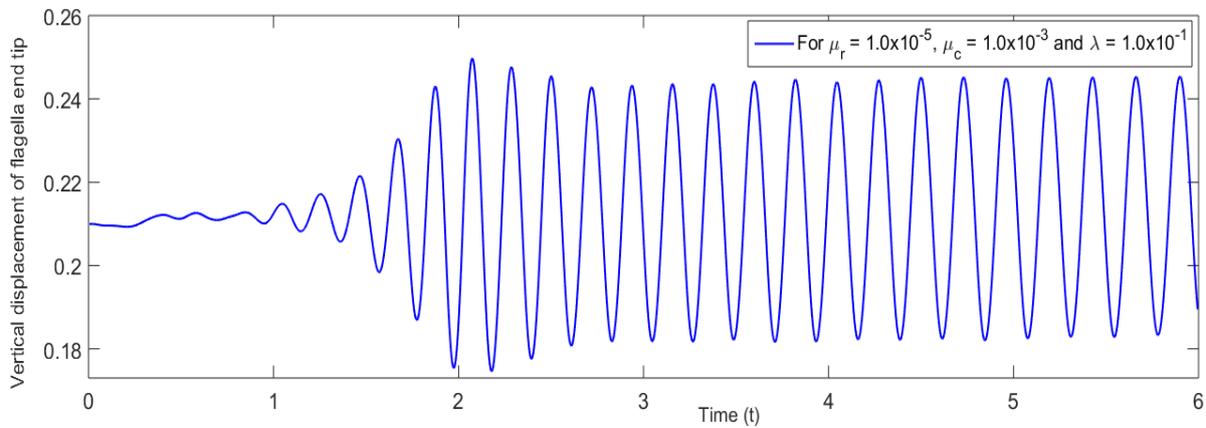

**Figure 4**: *Vertical displacement of flagella end tip (control point A) against microrotional viscosity $\mu_r$, microinnertia $\mu_c$ at different time. The frequency and amplitude is around $4.5 s^{-1}$ and $0.03$, respectively. The mesh has $2199$ and the time step is $0.005$*

Moreover, the structural material parameter $c_1$ plays a significant role and effects amplitude of oscillations as in Figure 5. It is found that the amplitude of oscillation decreases with increasing value of material parameter $c_1$ and larger mesh can make this relationship smooth as shown in Figure 6.

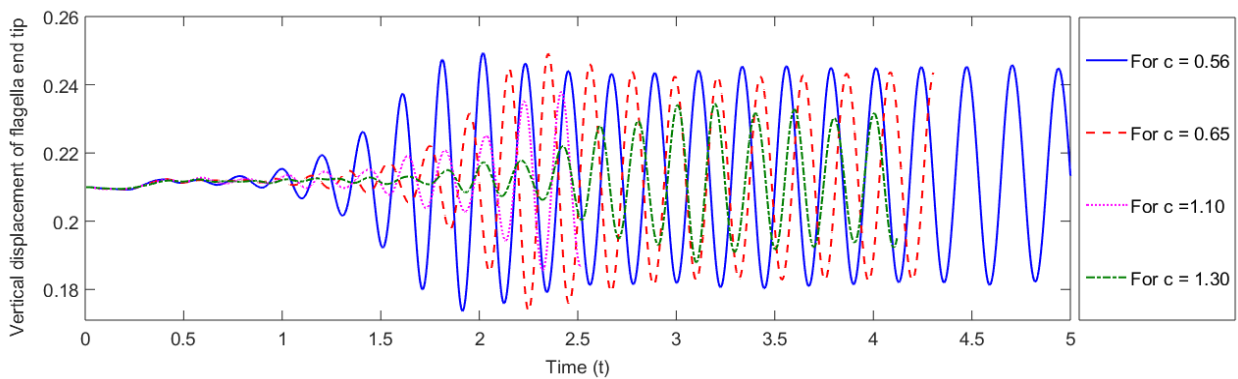

**Figure 5**: *Vertical displacement of flafella end tip (control point A) against different values of material parameter $c_1$.*

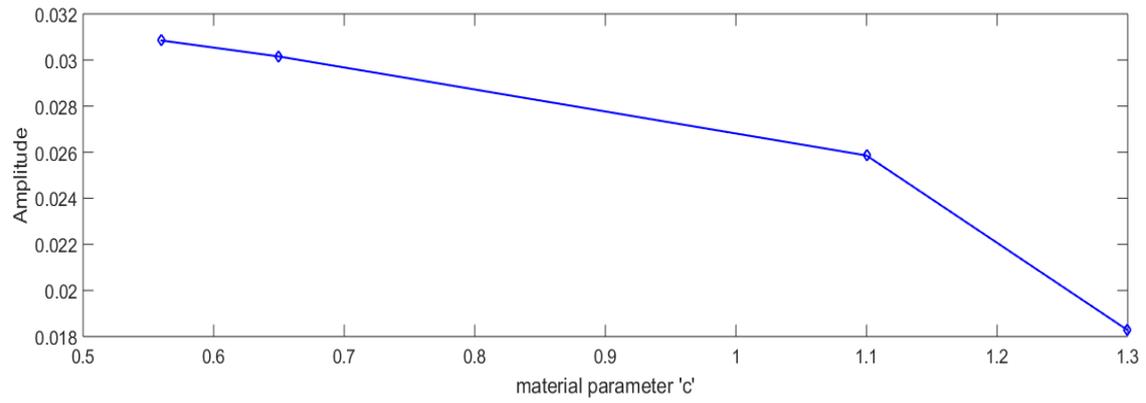

*Figure 6*: Change in amplitude with increasing value of material parameter $c_1$

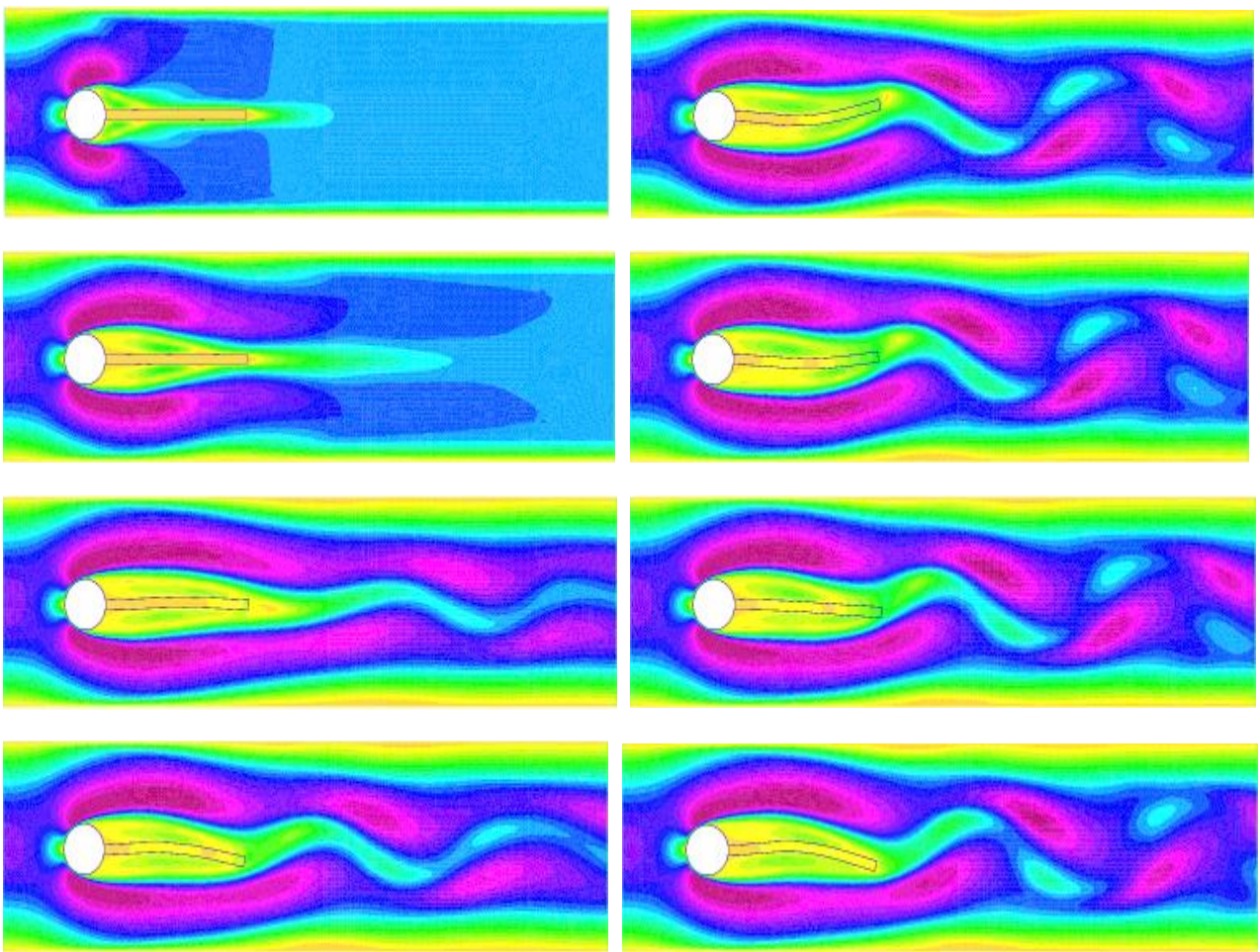

*Figure 7*: Velocity profile $V_h$ plots (*column* $1: t = 0.1, 0.3, 1.95, 2.3$ & *column* $1: t = 2.44, 2.46, 2.48, 2.52$ ).

Furthermore, the non-classical case of Cosserat fluid-structure interaction CFSI has been studied and obtained results are compared with the classical case of FSI. The classicial case of FSI governing model equations can be obtained by vanishing micro-rotational coeffinients. In non-classical FSI problems fluid particles exhibits micro-rotational behavior. The velocity compnents for classical and non-classical case are compared. It is observed that increasing values of micro-rotional viscosity develop more rotional effect on velocity components which validates the non-classical phenomena for CFSI problems. The graphical representation of the results in this case are displayed in Fiures 8-10.

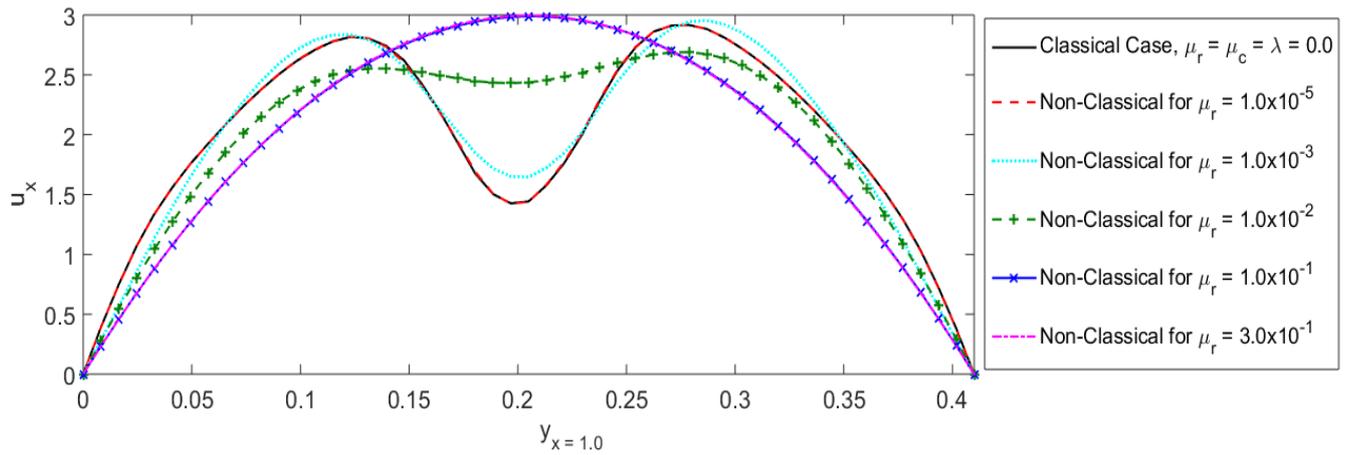

**Figure 8**: *Comparison of classical and non-classical case for $x-$components of velocity.*

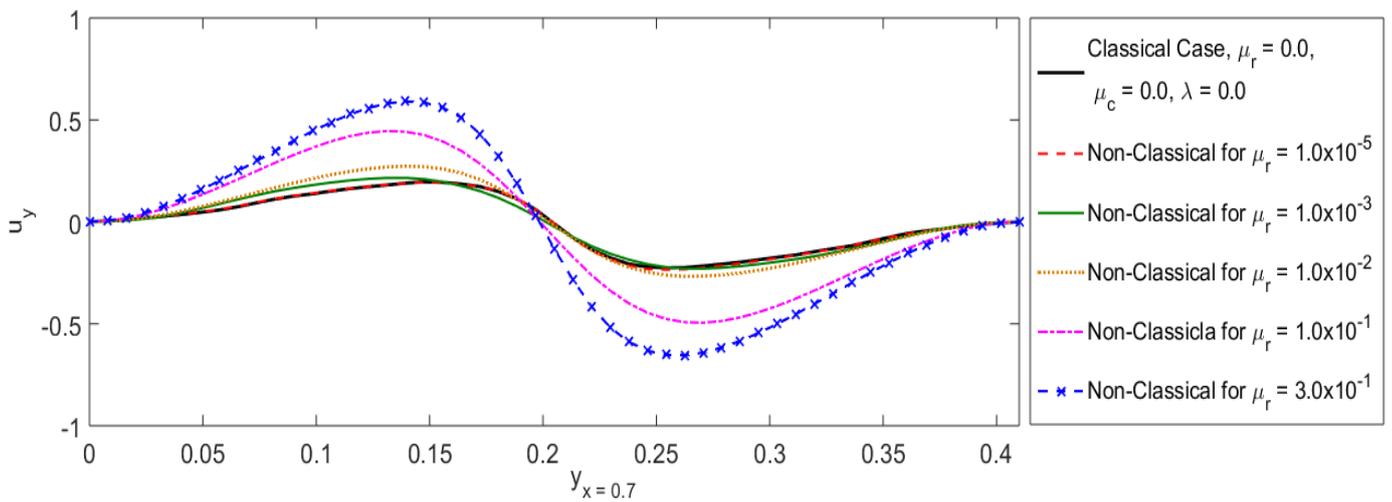

**Figure 9**: *Comparison of classical and non-classical case for $y-$ components of velocity.*

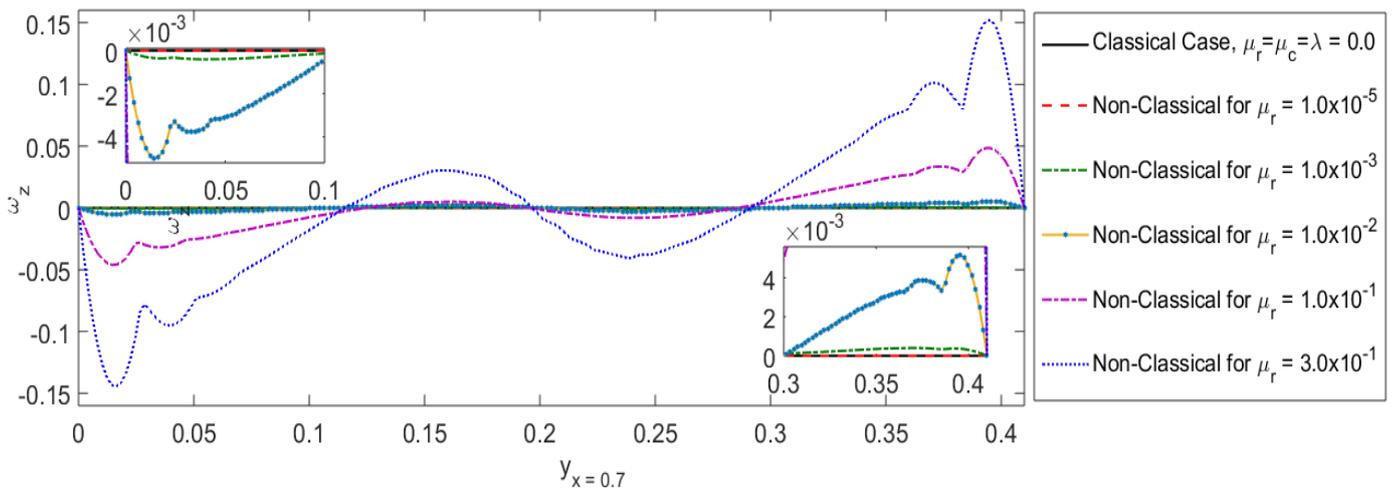

**Figure 10**: *Comparison of classical and non-classical case for $z-$components of micro-rotational velocity.*

# 7 CONCLUSION

In this research a monolithic Eulerian Cosserat fluid-structure interaction CFSI formulation has been presented in non-classical sense. The finite element method and semi-implicit scheme are used for discretizing space and time. The method has been implemented with a public domain software FreeFem ++. The graphical representation of results and color maps of the simulation at different time are also presented. The rotational effect at micro-structure level for CFSI has been studied and results are have been validated and compared with classical case of FSI. It is observed that amplitude of the oscillation decreases with increasing value of the material parameter $c_1$. The micro-rotational fluid parameters and their effect will be studied and analyzed in the second part of this paper.

**Author's Contribuitions:** Conceptualization, Muhammad Sabeel Khan; Writing, Muhammad Sabeel Khan and Nazim Hussain and Lisheng Liu; Methodology, Muhammad Sabeel Khan; Supervision, Muhammad Sabeel Khan.